\documentclass[5p]{elsarticle}
\usepackage{amsmath}
\usepackage{amsthm}
\usepackage{amssymb}
\usepackage{graphicx}
\usepackage{verbatim}
\usepackage{alltt}
\usepackage{color}

\newtheorem{thm}{Theorem}
\newtheorem{lem}[thm]{Lemma}
\newdefinition{rmk}{Remark}
\newproof{pf}{Proof}

\newcommand{\bfa}{\mathbf{a}}
\newcommand{\bfp}{\mathbf{p}}
\newcommand{\bfq}{\mathbf{q}}
\newcommand{\bfr}{\mathbf{r}}

\newcommand{\bfv}{\mathbf{v}}

\definecolor{red}{rgb}{0.4,0,0}
\definecolor{grn}{rgb}{0,0.4,0}
\definecolor{blu}{rgb}{0,0,0.6}

\journal{Computer Aided Geometric Design}

\begin{document}

\begin{frontmatter}

\title{A minimum-energy quadratic curve through three points and corresponding cubic Hermite spline}

\author[csu]{Steven Benoit\fnref{fn1}}
\ead{benoit@math.colostate.edu}
\address[csu]{Department of Mathematics, Colorado State University, Fort Collins, CO 80523 USA}
\fntext[fn1]{This research was supported by the National Science Foundation under grant GDE-0841259.}

\begin{abstract}
We demonstrate a method for exact determination of the quadratic curve of minimal energy and minimal curvature variation through three non-colinear points in the plane, including methods to determine the tangent vector and curvature at any point along the curve and an exact expression for the arc length of the curve between the first and last points.  We then extended this to a novel method of selecting tangent vectors for use in constructing Hermite splines to smoothly interpolate ordered sets of control points.  Our results are spline curves of lower energy than that of many popular spline implementations in most cases, which a series of examples demonstrate.
\end{abstract}

\end{frontmatter}

\section{Introduction}
Interpolating splines have long been used to provide smooth curves through discrete sets of points. The applications of such curves include interpolating data sets, computer-aided geometric design, typography, path planning for robotics and CNC machining, highway and railway design, and in the design of computer games \cite{Sc1973, Sp1995, Bo2001, Fa2002, Ma2005, PiGuRo2007, WaMe2005, Eb2007}.  Interpolating splines pass through each of an ordered set of defined control points, as opposed to estimating splines, which provide smooth curves that pass near their control points.

There are many forms of interpolating spline, some more apt to a specific application than others.  Two widely used measurements of the quality of an interpolating spline are the \emph{energy} (the elastic energy of the spline if it were an isomorphic flexible beam), and the \emph{curvature variation} (how much does the curvature change over the length of the spline).  Some general surveys, for example \cite{LeSe2009}, define additional criteria to judge interpolating spline quality, or \emph{fairness}, including
\begin{itemize}
  \item extensionality (does addin\textit{}g a new data point on the spline change the spline?),
  \item roundness (does the spline generate a circle given points lying on a circle?),
  \item monotone curvature (do the curvature extrema fall on the control points?), and
  \item locality vs. order (increasing order decreases ability to control spline locally).
\end{itemize}
Some examples of interpolating splines currently in use or under research include
\begin{itemize}
  \item The family of parametric cubic curves \cite{FoVaFeHu1996, Mo2006}, including Hermite curves (of which Catmull-Rom \cite{CaRo1974} and Kochanek-Bartels \cite{KoBa1984} splines are special cases) and B\'{e}zier curves
  \item Pythagorean-hodograph curves \cite{FaNe1995}
  \item Minimal-energy curves \cite{Ho1983}
  \item Minimal-variation curves \cite{Mo1992}
  \item Euler's spiral (or Conru spiral or clothoid curve) \cite{Le2008}
  \item Log-aesthetic curves \cite{YoFuSa2009}
\end{itemize}
Each of these have different strengths with respect to the criteria listed above.  However, the family of spline curve that is most commonly used in design work are the parametric curves, with cubics being the most common.

There are an infinite number of such curves that can interpolate any given data set.  In the case of Hermite curves, the choice of tangent vectors at each point selects from this family of curves, and the quality of the resulting interpolating curve is based strongly on the method of choosing these tangent vectors.  Before we discuss our method for selecting these tangent vectors, we briefly explore parametric quadratic curves.

\section{Quadratic parametric curves}

By parametric quadratic curves, we mean parametric curves represented by functions of the form
\begin{equation} \label{parametric-curve}
  \bfr(t) = \bfa_1 t^2 + \bfa_2 t + \bfa_3 \, .
\end{equation}
The curvature $\kappa$, elastic energy $E$ and curvature variation $V$ of a parametric curve are given by
\begin{equation} \label{kappa_def}
  \kappa = \frac{\dot{x} \ddot{y} - \dot{y} \ddot{x}}
                {(\dot{x}^2 + \dot{y}^2)^{3/2}} \, ,
\end{equation}
\begin{equation} \label{energy}
  E = \int \kappa(t)^2 \, \mbox{d}t \, ,
\end{equation}
\begin{equation} \label{curv_variation}
  V = \int \dot{\kappa}(t)^2 \, \mbox{d}t \, ,
\end{equation}
where dots indicate derivatives with respect to the curve parameter.  The first and second derivatives of (\ref{parametric-curve}) are
\begin{equation} \label{r_derivatives}
  \dot{\bfr}(t) = 2 \bfa_1 t + \bfa_2
  \qquad \mbox{and} \qquad
  \ddot{\bfr}(t) = 2 \bfa_1 \, ,
\end{equation}
and using $\bfa_i = (x_i, y_i), i \in \{1,2,3\}$, we have
\begin{equation} \label{kappa}
  \kappa(t) = 2 \frac{x_2 y_1 - x_1 y_2}
                     {\left[ (2 x_1 t + x_2)^2 + (2 y_1 t + y_2)^2 \right]^{3/2}} \, .
\end{equation}
Substituting (\ref{kappa}) into (\ref{energy}),
\[
  E = 4 (x_2 y_1 - x_1 y_2)^2
      \int \bigl[ (2 x_1 t + x_2)^2 + (2 y_1 t + y_2)^2 \bigr]^{-3} \mbox{d}t \, .
\]
Evaluating this integral over the whole real line yields
\[
  E = \frac{3 \pi}{4} \frac{(y_1^2 + x_1^2)^2} {|x_1 y_2 - x_2 y_1|^3} \, ,
\]
or in a coordinate-free form,
\begin{equation}
  E = \frac{3 \pi}{4} \frac{|\bfa_1|^4}{|\bfa_1 \times \bfa_2|^3}.
\end{equation}

Then, differentiating (\ref{kappa_def}) and using the fact that third derivatives vanish for a quadratic curve,
\begin{equation}
  \dot{\kappa} =
    -3 \frac{(\dot{x} \ddot{y} - \dot{y} \ddot{x})
             (\dot{x} \ddot{x} + \dot{y} \ddot{y})}
            {(\dot{x}^2 + \dot{y}^2)^{5/2}} \, .
\end{equation}
Applying (\ref{r_derivatives}),
\[
  \dot{\kappa} =
     12 \frac{\left( x_1 y_2 - x_2 y_1 \right)
              \left( 2 \left(x_1^2 +  y_1^2\right) t + x_1 x_2 + y_1 y_2 \right)}
             {\left[ \left(2 x_1 t + x_2 \right)^2 + \left(2 y_1 t + y_2 \right)^2 \right]^{5/2}} \, ,
\]
and so the curvature variation is
\begin{align*}
  V = 144 &\left( x_1 y_2 - x_2 y_1 \right)^2 \\
  &
      \int \frac{\left[2 \left(x_1^2 +  y_1^2\right) t + x_1 x_2 + y_1 y_2 \right]^2}
                {\left[ \left(2 x_1 t + x_2 \right)^2 + \left(2 y_1 t + y_2 \right)^2 \right]^5}
      \, \mbox{d}t \, .
\end{align*}
Evaluating this over the real line gives the result
\begin{align*} \label{curv-var}
  V =& \frac{45 \pi (x_1^2 + y_1^2)^2}
        {16 |x_1 y_2 - x_2 y_1|^7}
             \bigl[ (x_1 + y_1^2)^2 |x_1 y_2 - x_2 y_1|^2 \\
& \qquad \qquad \qquad \qquad
                   + 7 x_1^2 (x_1 - 1)^2 (x_1^2 + y_1 y_2)^2
             \bigr]\, .
\end{align*}
However, we are more interested in minimizing $V$ rather than obtaining its exact value, and so we introduce the following useful Lemma,
\begin{lem}
A quadratic curve of least energy also has least curvature variation.
\begin{pf}
Consider a quadratic curve defined parametrically by $y = a t^2 + b t + c$, with tangent vector given by $\dot{y} = 2 a t + b$.
The curvature at a point $t$ is given by
\[
  \kappa = 2 a \left[ 1 + (2 a t + b)^2 \right]^{-3/2} \, ,
\]
with rate of change given by
\[
  \dot{\kappa} = -12 a^2 (2 a t + b) \left[ 1 + (2 a t + b)^2 \right]^{-5/2} \, .
\]
The total energy $E$ and curvature variation $V$ of the curve are given by
\begin{equation} \label{eq:e}
  E = 4 a^2 \int_{-\infty}^\infty \left[ 1 + (2 a t + b)^2 \right]^{-3} \, \mbox{d}t \, ,
\end{equation}
and
\begin{equation} \label{eq:v}
  V = 144 a^4 \int_{-\infty}^\infty (2 a t + b)^2 \left[ 1 + (2 a t + b)^2 \right]^{-5} \, \mbox{d}t \, .
\end{equation}
Integrating (\ref{eq:e}) and (\ref{eq:v}),
we obtain
\[
  E = \frac{3 \pi}{4} |a|   \, ,
  \qquad \mbox{ and} \qquad
  V = \frac{45 \pi}{16} |a| \, .
\]
Therefore, the quadratic that minimizes $|a|$ will minimize both energy and curvature variation.
\end{pf}
\end{lem}

\section{Quadratic curve passing through three points}

Suppose we have three non-colinear points $\bfp_1, \bfp_2, \bfp_3 \in \mathbb{R}^2$.  As pointed out in \cite{LaSc1991}, these points can be interpolated by a parabola with $\bfr(0) = \bfp_1$, $\bfr(1) = \bfp_3$, and $\bfr(0.5) = \bfp_2$ (what Lachance and Schwartz call the ``Cinci Parabola'').  However, there is no reason to require that $r(t)$ pass $\bfp_2$ at $t=0.5$.  Rather, we wish to find the quadratic curve $\bfr(t)$ of least energy such that $\bfr(0) = \bfp_1$, $\bfr(1) = \bfp_3$, and $\bfr(T) = \bfp_2$ where $0 < T < 1$. Then,
\begin{align*}
  \bfp_1 =& \bfr(0) = \bfa_3 \, ,\\
  \bfp_2 =& \bfr(T) = \bfa_1 T^2 + \bfa_2 T + \bfa_3 \, , \quad \mbox{and} \\
  \bfp_3 =& \bfr(1) = \bfa_1 + \bfa_2 + \bfa_3 \, .
\end{align*}
Solving for $\bfa_1$, $\bfa_2$, and $\bfa_3$,
\begin{equation} \label{a_in_terms_of_t}
\begin{split}
  \bfa_1 =& \frac{\bfp_2 - \bfp_1 - (\bfp_3 - \bfp_1) T}{T^2 - T} \, , \\
  \bfa_2 =& \bfp_3 - \bfp_1 - \frac{\bfp_2 - \bfp_1 - (\bfp_3 - \bfp_1) T}{T^2 - T} \, , \quad \mbox{and} \\
  \bfa_3 =& \bfp_1 \, .
\end{split}
\end{equation}
We transform coordinates so $\widetilde{\bfp}_1$ lies at the origin, $\widetilde{\bfp}_i = \bfp_i - \bfp_1$, $i = 1, 2, 3,$
and then scale so $|\widehat{\bfp}_3 - \widehat{\bfp}_1| = 1$, using
\[
  \widehat{\bfp}_i = \frac{\widetilde{\bfp}_i}{|\widetilde{\bfp}_3 - \widetilde{\bfp}_1|}
                         = \frac{\bfp_i - \bfp_1}{|\bfp_3 - \bfp_1|}
  \, , \qquad  i = 1, 2, 3.
\]
We now use $\widehat{\bfp}_3$ to generate rotation matrix $A$ that will rotate that point to $(1,0)$,
\[
  A = \left|
    \begin{array}{c c}
       {\mbox{$\hat{p}_3$}}_x & {\mbox{$\hat{p}_3$}}_y \\
      -{\mbox{$\hat{p}_3$}}_y & {\mbox{$\hat{p}_3$}}_x \\
    \end{array}
    \right| \, .
\]
Applying $A$ to each point generates working points $\bfq_i = A \widehat{\bfp}_i$.  In this new coordinate frame, $\bfq_1 = 0$, $\bfq_3 = (1,0)$, and ${q_2}_y \ne 0$.  Equations (\ref{a_in_terms_of_t}) are equally valid in this new frame, in which case they simplify somewhat,
\begin{equation} \label{a_vectors}
\begin{split}
  \bfa_1 =& \frac{\bfq_2 - \bfq_3 T}{T^2 - T} \, , \\
  \bfa_2 =& \bfq_3 - \frac{\bfq_2 - \bfq_3 T}{T^2 - T} \, , \quad \mbox{and} \\
  \bfa_3 =& \mathbf{0} \, .
\end{split}
\end{equation}
\subsection{Minimum energy quadratic}

The energy, in terms of $T$, is
\[
  E = \frac{3 \pi \alpha}{4 \left| \bfq_2 \times \bfq_3 \right|^3}
      \frac{|\bfq_2 - \bfq_3 T|^4} {T^2 - T} \, .
\]
To find the value of $T$ for which the curve has minimal energy, we set,
\[
  \frac{\partial}{\partial T}
  \frac{\left[ (\bfq_2 - \bfq_3 T) \cdot (\bfq_2 - \bfq_3 T) \right]^2}
       {T^2 - T} = 0 \, ,
\]
Solving for $T$, recalling that $|\bfq_3| = 1$, generates a cubic equation in $T$,
\begin{equation} \label{cubic_in_t}
     T^3
   - \frac{3}{2} T^2
   + ({q_2}_x - |\bfq_2|^2) T
   + \frac{1}{2} |\bfq_2|^2
   = 0 \, .
\end{equation}
When $T=0$, the left-hand side of (\ref{cubic_in_t}) is positive definite, and when $T=1$, the left-hand side reduces to
\[
   -\frac{1}{2} \left[ (1 - {q_2}_x)^2 + {{q_2}_y}^2 \right] \, ,
\]
which is negative definite since ${q_2}_y \ne 0$.  Therefore, by continuity of (\ref{cubic_in_t}), we can be assured of having a root in the range $0 < T < 1$.  Moreover, since (\ref{cubic_in_t}) goes to $\infty$ as $T \to \infty$ and goes to $-\infty$ as $T \to -\infty$, (\ref{cubic_in_t}) has three real roots, and we seek the middle root of the three.

The roots of the cubic can be found using the cubic formula,
\begin{equation} \label{roots_solution}
\begin{split}
  T_1 =& \frac{1}{2} + \frac{\mu + \zeta}{2}, \\
  T_2 =& \frac{1}{2} - \frac{\mu + \zeta - i \sqrt{3} (\mu - \zeta)}{4} \, , \quad \text{and} \\
  T_3 =& \frac{1}{2} - \frac{\mu + \zeta + i \sqrt{3} (\mu - \zeta)}{4} \, ,
\end{split}
\end{equation}
where
\begin{align*}
  &\mu    = \sqrt[3] {\beta + \sqrt{\gamma + \beta^2}} \, , \quad
   \zeta  =  \sqrt[3] {\beta - \sqrt{\gamma + \beta^2}} \, , \\
  &\beta  =  1 - 2 {q_2}_x \, , \quad \text{and} \quad
   \gamma =  \frac{\left[ 4 ({q_2}_x - |\bfq_2|^2) - 3 \right]^3}{27} \, .
\end{align*}
For example, suppose $\bfq_2 = \left( \frac{1}{2}, 1 \right)$, in which case we would expect $T = \frac{1}{2}$ by symmetry.
In this case, $\beta=0$ and $\gamma = -8$, giving $\mu = i \sqrt{2}$ and $\zeta = -i \sqrt{2}$, so the roots are
\[
  T_1 = \frac{1}{2}, \quad
  T_2 = \frac{1}{2} + \frac{\sqrt{6}}{2}, \quad \text{and} \quad
  T_3 = \frac{1}{2} - \frac{\sqrt{6}}{2} \, ,
\]
and indeed, we obtain three real roots, one of which lies in $(0,1)$, that root falling at $T = \frac{1}{2}$ as expected.

Once $T$ (the root between $0$ and $1$) has been identified, we use (\ref{a_in_terms_of_t}) to recover the $\bfa_i$ of the desired minimum energy curve, then (\ref{parametric-curve}) to generate the curve.  Reversing the coordinate transformations is not necessary since the value $T$ is independent of coordinate system.  A summary of the algorithm presented here, in more convenient pseudocode form, is included in the Appendix.

\section{Tangent vector at point $\bfp_2$}

To construct a Hermite spline through a sequence of points in which $(\bfp_1, \bfp_2, \bfp_3)$ is a subsequence, we construct tangent vector $\dot{\bfr}(T)$ at $\bfp_2$ to the minimal-energy quadratic found above.  Having solved for $T$, we have
\[
  \dot{\bfr}(T) = (2T - 1) \left[ \frac{\bfp_2 - \bfp_1 - (\bfp_3 - \bfp_1) T}{T^2 - T} \right] + \bfp_3 - \bfp_1
\]

\begin{figure}[!ht]
 \centering
 \begin{tabular}{c c}
  \includegraphics[width=0.17\textwidth]{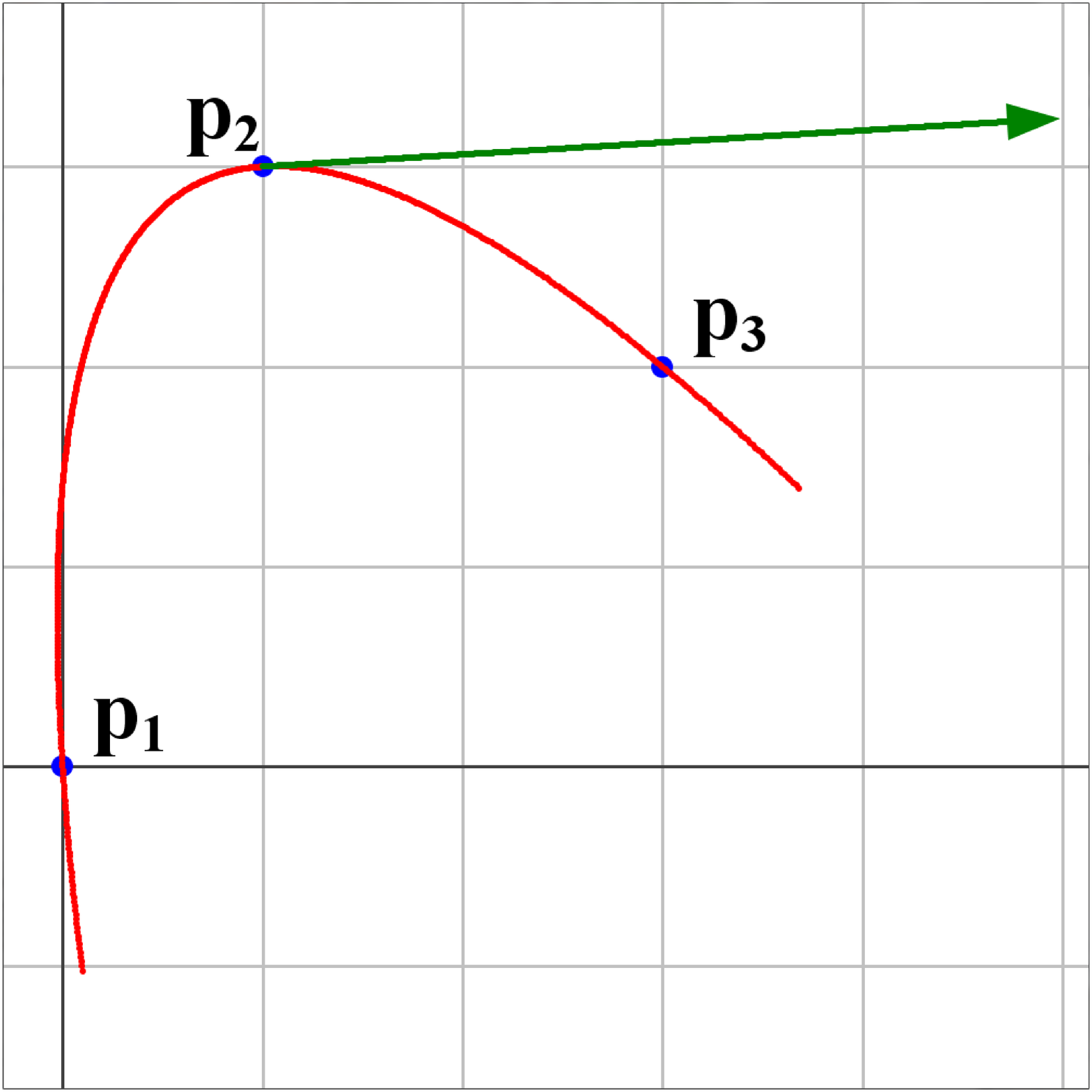} &
  \includegraphics[width=0.17\textwidth]{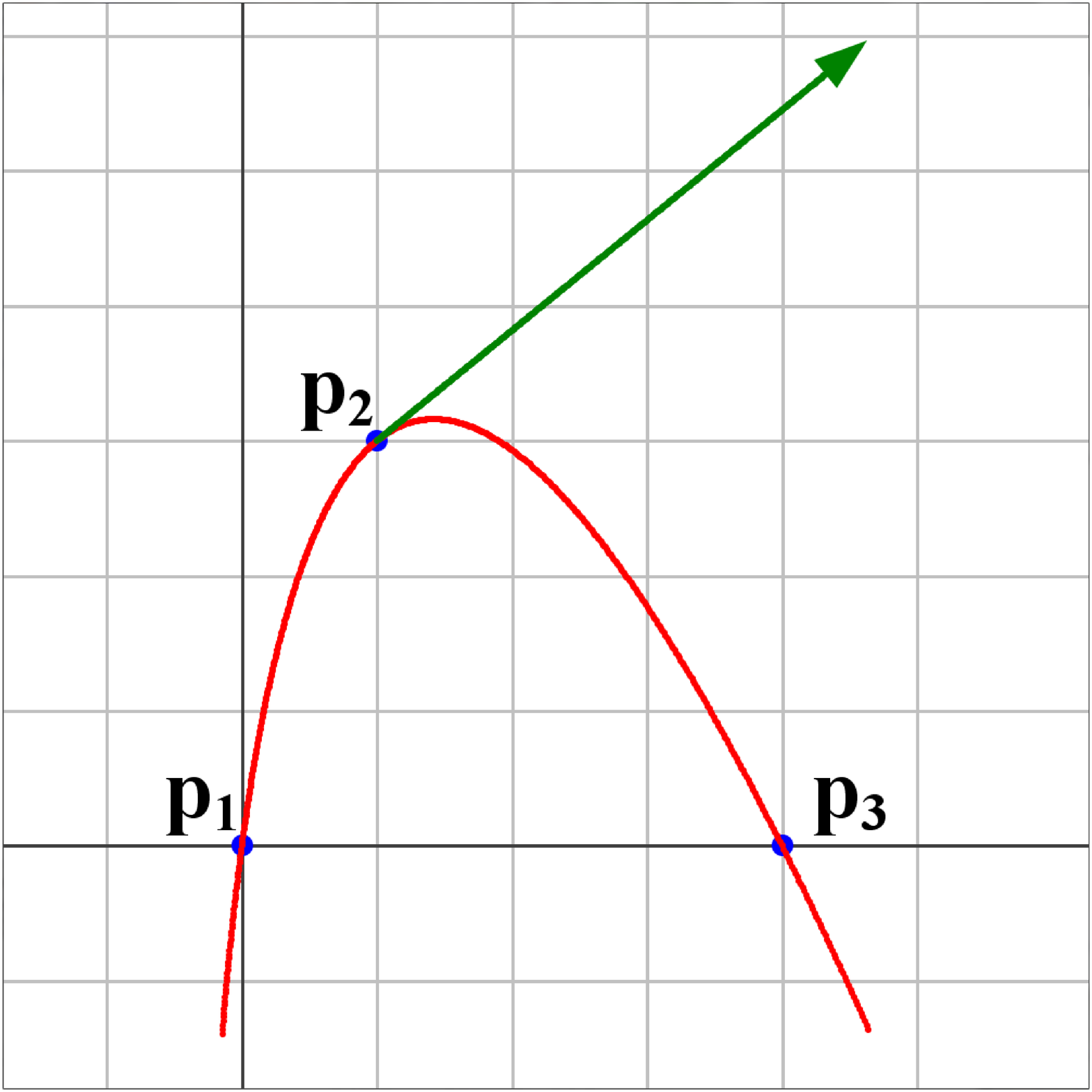} \\
  \includegraphics[width=0.17\textwidth]{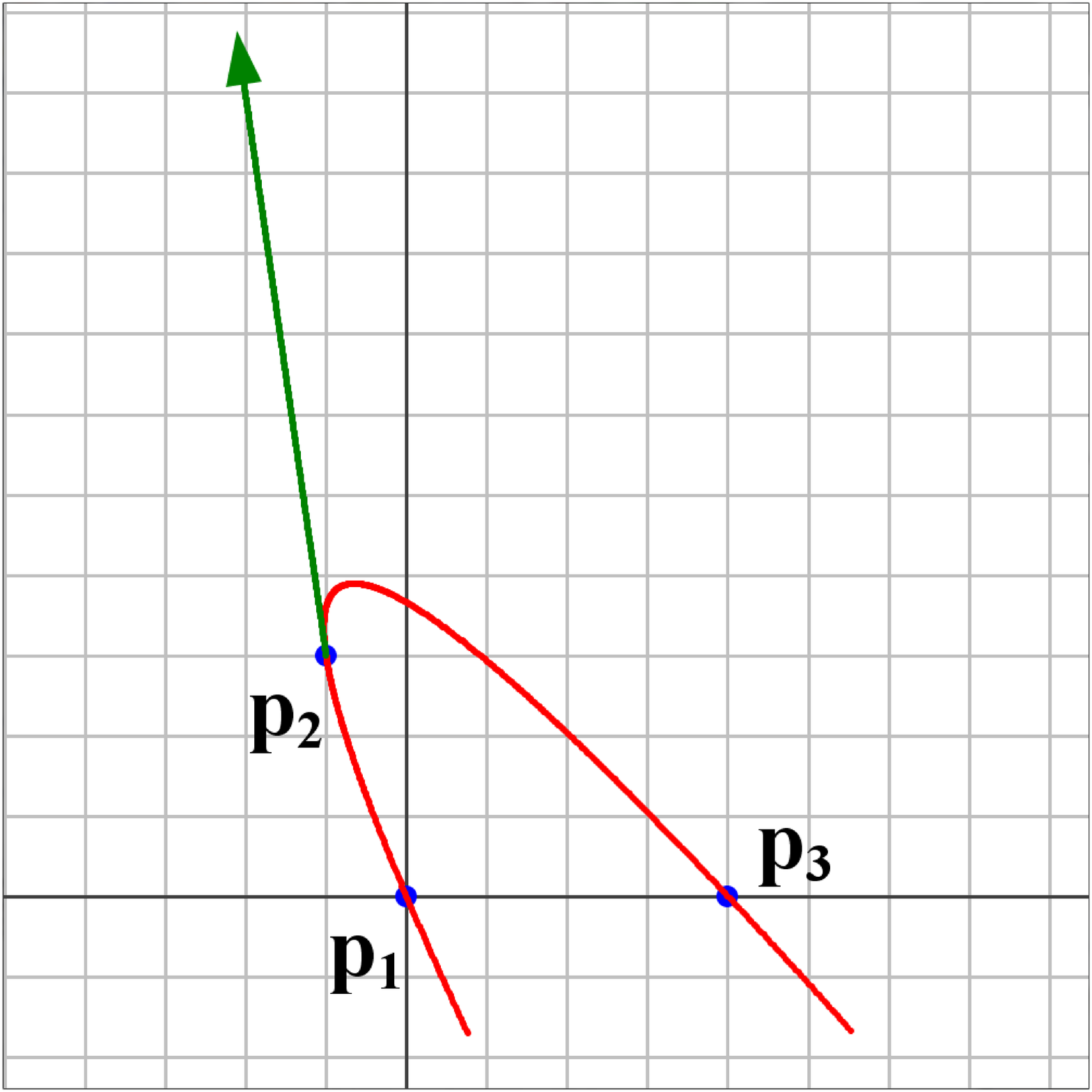} &
  \includegraphics[width=0.17\textwidth]{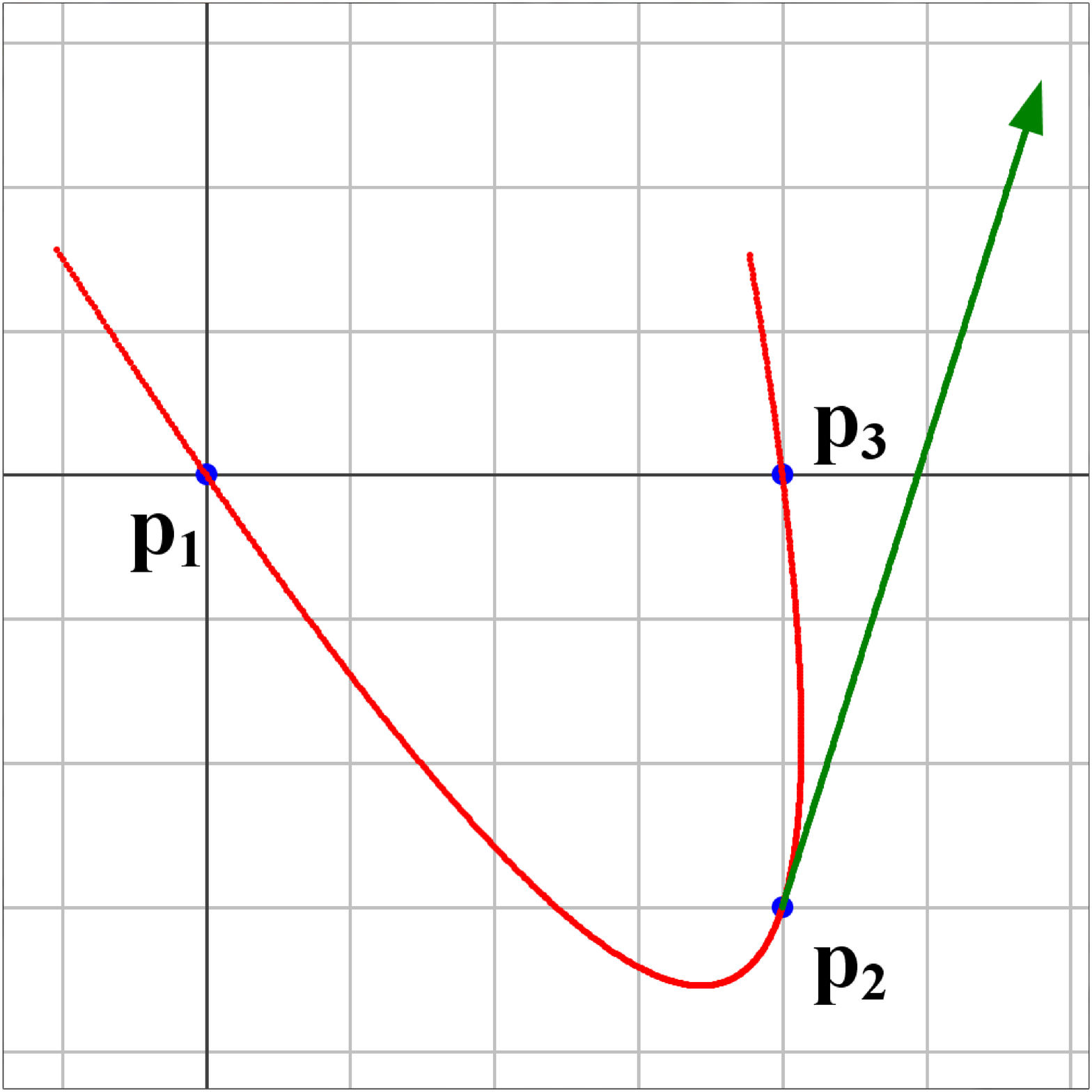}
 \end{tabular}
 \caption{Examples of ordered three-point sets and the corresponding minimum-energy quadratic.  The tangent vector (derivative of the curve with respect to its parameter) at the center point is also shown.}
 \label{fig:examples}
\end{figure}

For example, if $\mathbf{p}_1=(0,0)$, $\mathbf{p}_2=(\frac{1}{2},1)$, and $\mathbf{p}_1=(1,0)$, then $\dot{\bfr}(T) = (1,0)$.  Some examples of minimum-energy quadratic functions and the corresponding tangent vector are shown in Figure \ref{fig:examples}.

\section{Arc length between $\bfp_1$ and $\bfp_3$}

It may also be useful to know the arc length $l$ of the curve between points $\bfp_1$ and $\bfp_3$.  The arc length is found by integrating the norm of the tangent vector over the interval.

\[
\begin{split}
  l =& \int_0^1 |2 t \bfa_1 + \bfa_2| \, \mbox{d}t \\
    =& \int_0^1 \Biggl|
            2 t \left( \frac{\bfp_2 - \bfp_1 - (\bfp_3 - \bfp_1) T}{T^2 - T} \right) \\
& \qquad \qquad
          + \bfp_3 - \bfp_1 - \frac{\bfp_2 - \bfp_1 - (\bfp_3 - \bfp_1) T}{T^2 - T}
        \Biggr| \, \mbox{d}t \, .
\end{split}
\]
We define $\mathbf{s}_3 = \bfp_3 - \bfp_1$, and $\mathbf{s}_2 = \bfp_2 - \bfp_1$, so
\[
\begin{split}
  l =& \int_0^1 \left|
         (2 t - 1) \left( \frac{T \mathbf{s}_3 - \mathbf{s}_2}{T - T^2} \right) + \mathbf{s}_3
     \right| \, \mbox{d}t \\
    =& \int_0^1 \left| \frac{1}{T - T^2} \left[
         (2 t - 1) \left( T \mathbf{s}_3 - \mathbf{s}_2 \right) + (T - T^2) \mathbf{s}_3 \right]
     \right| \, \mbox{d}t \, .
\end{split}
\]
Then we note that for $0 < T < 1$, $T - T^2 > 0$, so
\[
  l
    = \frac{1}{T - T^2} \int_0^1 \left|
           2 t \left( T \mathbf{s}_3 - \mathbf{s}_2 \right)
         - (T^2 \mathbf{s}_3 - \mathbf{s}_2)
      \right| \, \mbox{d}t \, .
\]
Defining $\bfr_1 = T \mathbf{s}_3 - \mathbf{s}_2$ and $\bfr_2 = T^2 \mathbf{s}_3 - \mathbf{s}_2$, and $\theta$ as the angle between $\bfr_1$ and $\bfr_2$,
\[
  l
    = \frac{1}{T - T^2} \int_0^1
      \sqrt{ 4 t^2 |\bfr_1|^2 - 4 t |\bfr_1| |\bfr_2| \cos \theta + |\bfr_2|^2} \, \mbox{d}t\, . 
\]
Evaluating this integral, we obtain
\[
\begin{split}
l =
\frac{1}{4 |\bfr_1|(T - T^2)}
   \biggl( &
       |\bfr_2|^2 \cos \theta
     + \left( 2 |\bfr_1| - |\bfr_2| \cos \theta \right)
       \rho + \\
& 
     |\bfr_2|^2 \sin^2 \theta
            \log \frac{  2 |\bfr_1|
                       - |\bfr_2| \cos \theta
                       + \rho}
                      {|\bfr_2| (1 - \cos \theta)}
   \biggr) \, ,
\end{split}
\]
where $\rho = \sqrt{4 |\bfr_1|^2 - 4 |\bfr_1| |\bfr_2| \cos \theta + |\bfr_2|^2}$.

\section{A cubic Hermite spline implementation}

A Hermite cubic curve between any two consecutive points in an ordered point set is based on the locations of the points and on a tangent vector assigned at each point.  The choice of tangent vectors at each vertex strongly affects the resulting spline.  We now use the minimum-energy quadratic derived above to generate these tangent vectors, and compare the results to other methods.  Note that the use of a single tangent vector at a control point for both adjoining spline segments assures $G^1$ continuity of the resulting curve.

\begin{figure}[!ht]
 \centering
 \begin{tabular}{c c}
  \includegraphics[width=0.17\textwidth]{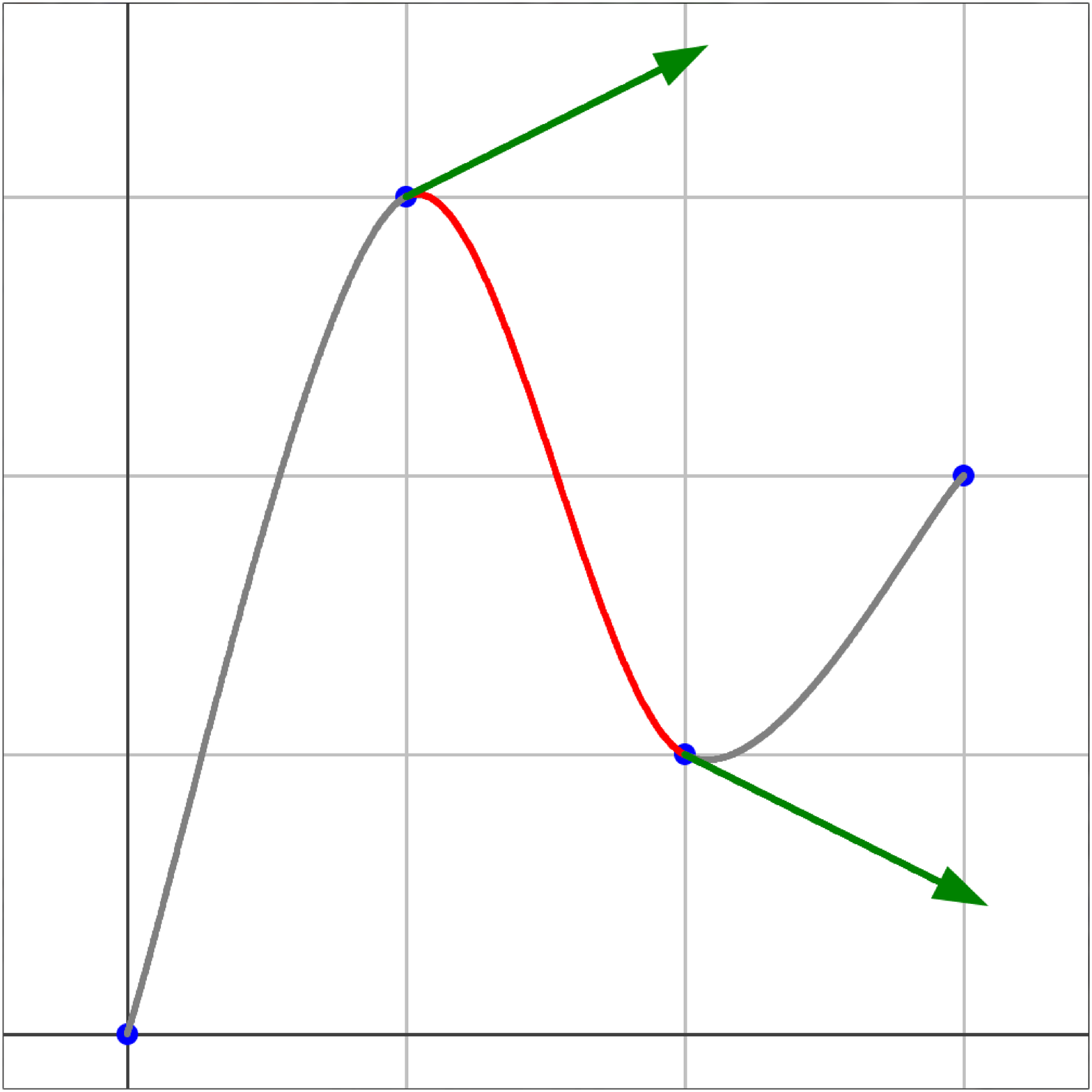} &
  \includegraphics[width=0.17\textwidth]{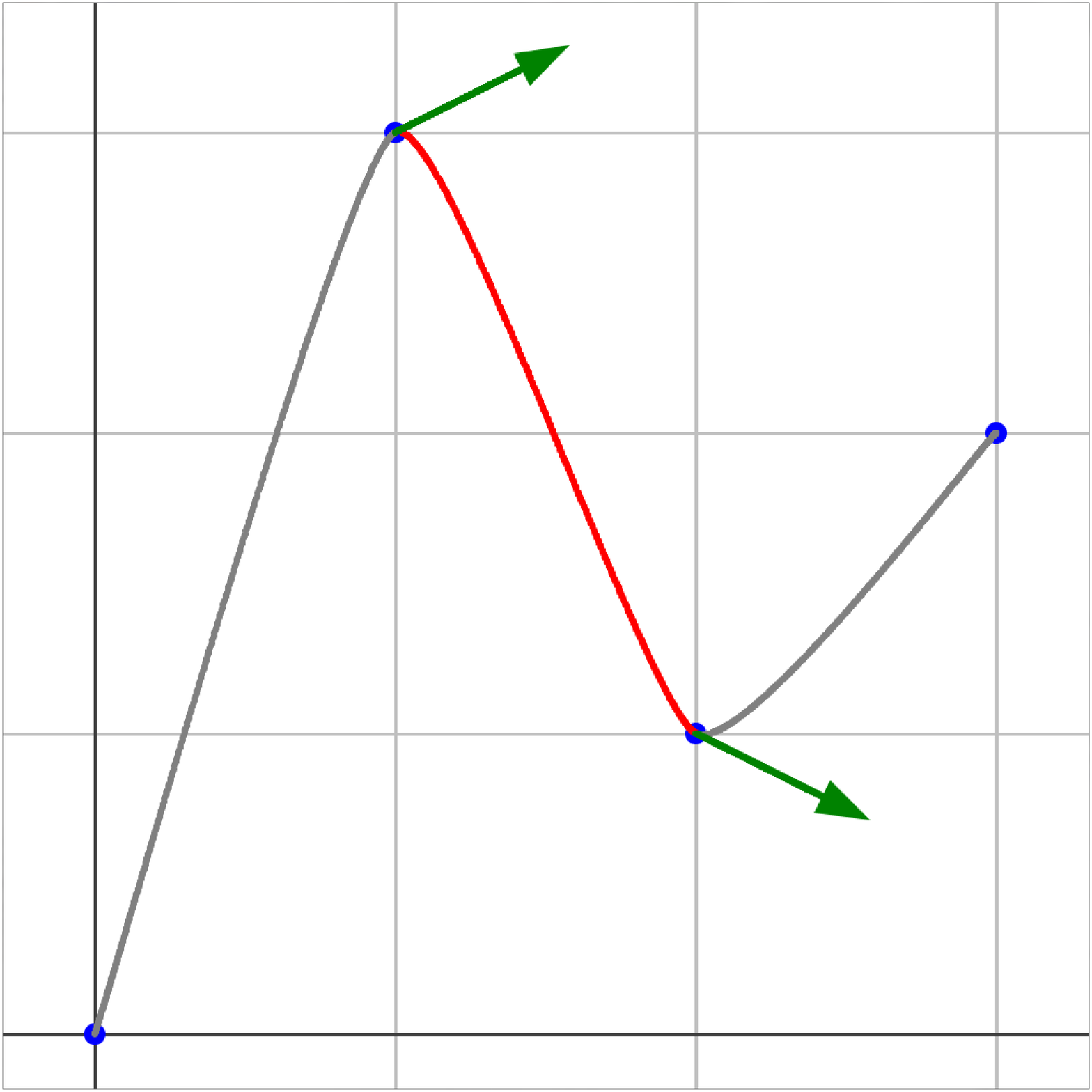} \\
  \footnotesize{Catmull-Rom} & 
  \footnotesize{Cardinal} \\
  \includegraphics[width=0.17\textwidth]{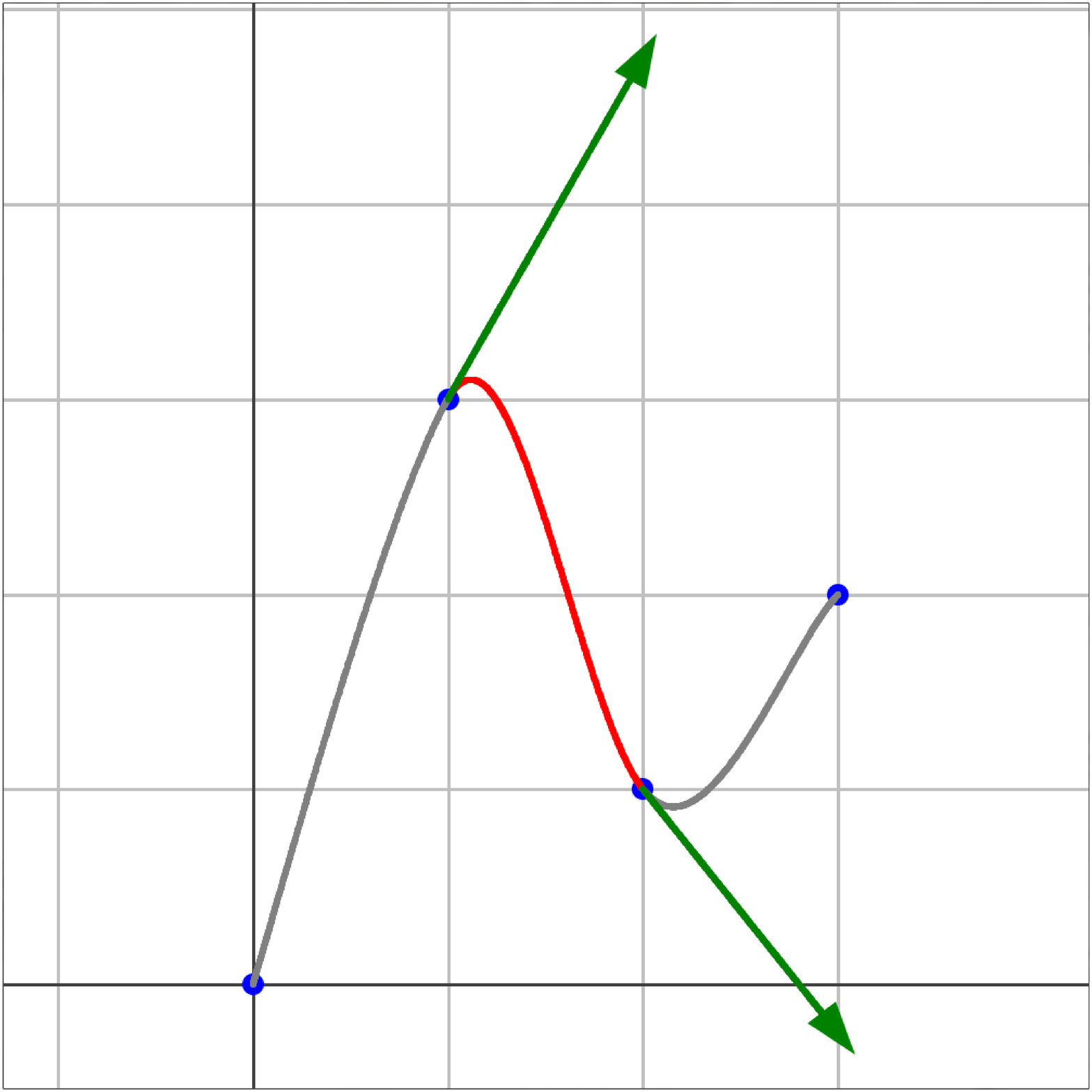} &
  \includegraphics[width=0.17\textwidth]{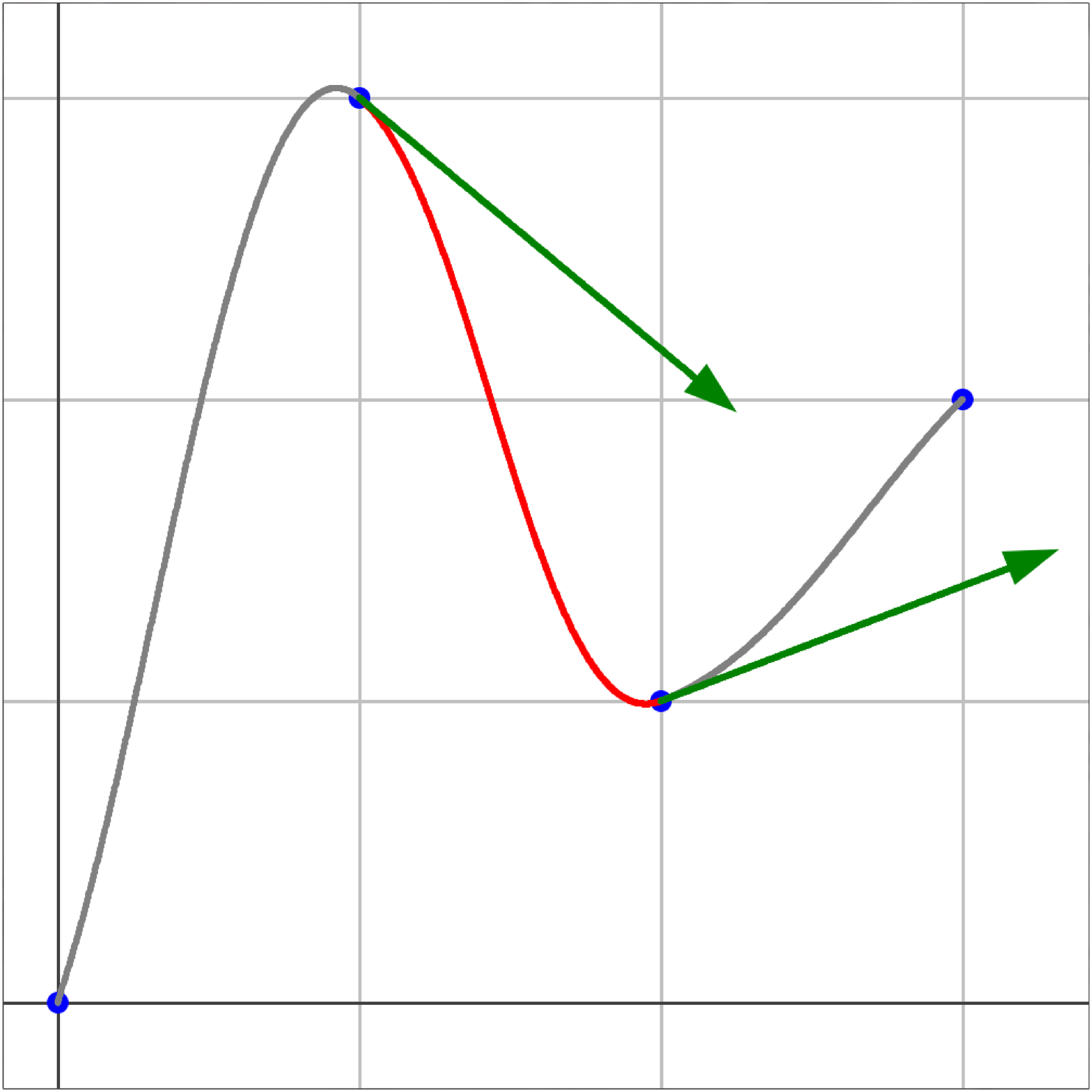} \\
  \footnotesize{Kochanek-Bartels} &
  \footnotesize{This paper's method}\\
 \end{tabular}
 \caption{Comparison of methods of choosing tangent vectors for cubic interpolating splines.  The cardinal spline shown uses a tension parameter of 0.5.  The Kochanek-Bartels spline uses a tension and continuity of 0 and bias of 0.5.}
 \label{fig:compare-splines}
\end{figure}

Figure \ref{fig:compare-splines} shows splines that result for several different methods of choosing of tangent vectors.  If control points $\{ \bfp_1, \bfp_2, ..., \bfp_n\}$ occur at parameter values $\{ t_1, t_2, ..., t_n\}$ along the spline, the methods of selecting tangent vectors $\{ \bfv_2, \bfv_3, ..., \bfv_{n-1}\}$ shown in the Figure are listed here.

\subsection{Catmull-Rom spline:}
\[
  \bfv_i = \frac{\bfp_{i+1} - \bfp_{i-1}}{t_{i+1} - t_{i-1}} \, ,
\]

\subsection{Cardinal spline:}
\[
  \bfv_i = (1-\tau) \frac{\bfp_{i+1} - \bfp_{i-1}}{t_{i+1} - t_{i-1}}
\]
where $\tau$ is a tension parameter,

\subsection{Kochabek-Bartels spline:}
\[
\begin{split}
  \bfv_i =& \frac{(1-\tau)(1+\beta)(1+\gamma)}{2} \left( \bfp_{i} - \bfp_{i-1} \right) \\
          & + \frac{(1-\tau)(1-\beta)(1-\gamma)}{2} \left( \bfp_{i+1} - \bfp_{i} \right)
\end{split}
\]
where $\tau$ is a tension parameter, $\beta$ is a bias parameter, and $\gamma$ is a continuity parameter.

\subsection{A spline based on minimum-energy quadratics:}
\[
\begin{split}
  \bfv_i = \frac{1}{t_{i+1} - t_{i-1}}
                 \biggl\{&
                    (2T-1) \left[
                              \frac{
                                       \bfp_i
                                     - \bfp_{i-1}
                                     - (\bfp_{i+1} - \bfp_{i-1}) T
                                   }
                                   {T^2 - T}
                           \right] \\
&
                    + \bfp_{i+1} - \bfp_{i - 1}
                 \biggr\}
\end{split}
\]
Where $T$ is the root from (\ref{roots_solution}) that lies in $0 < T < 1$.

\begin{figure}[!ht]
 \centering
 \begin{tabular}{c c}
 \includegraphics[width=0.17\textwidth]{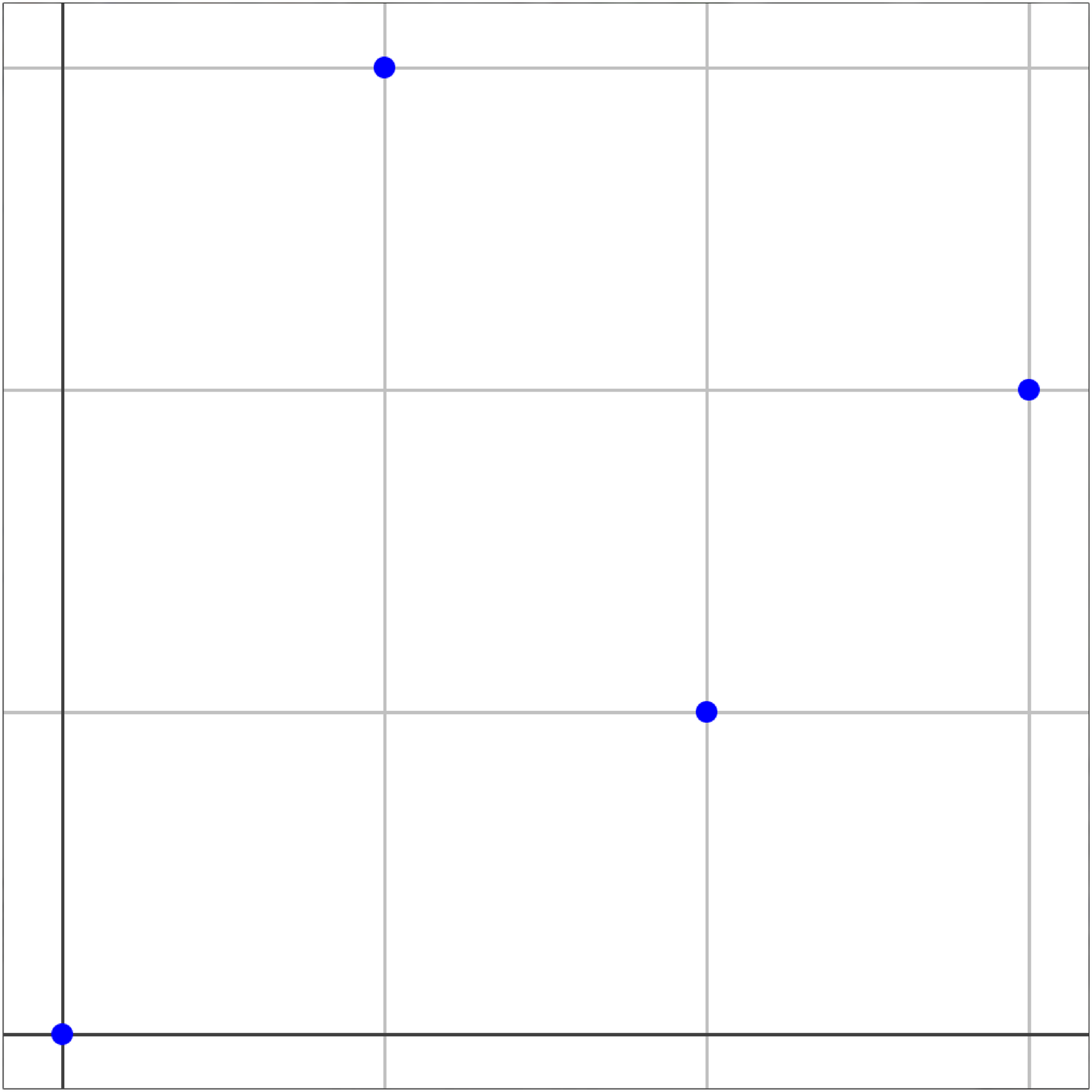} &
 \includegraphics[width=0.17\textwidth]{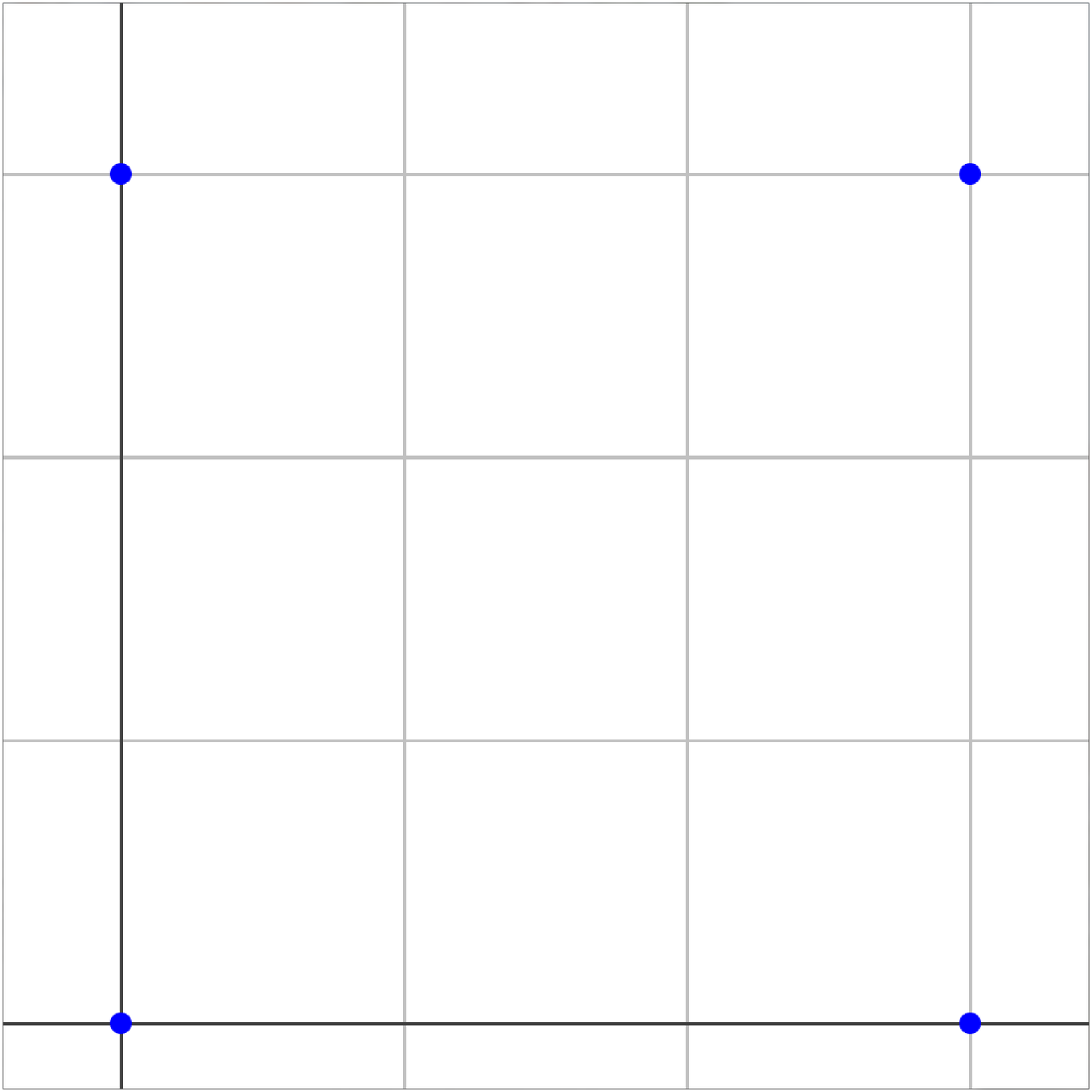} \\
 Point Set 1 & Point Set 2  \\
 \footnotesize{(0,0), (1,3), (2, 1), (3, 2)} &
 \footnotesize{(0,0), (0,3), (3, 3), (3, 0)} \\
 \includegraphics[width=0.17\textwidth]{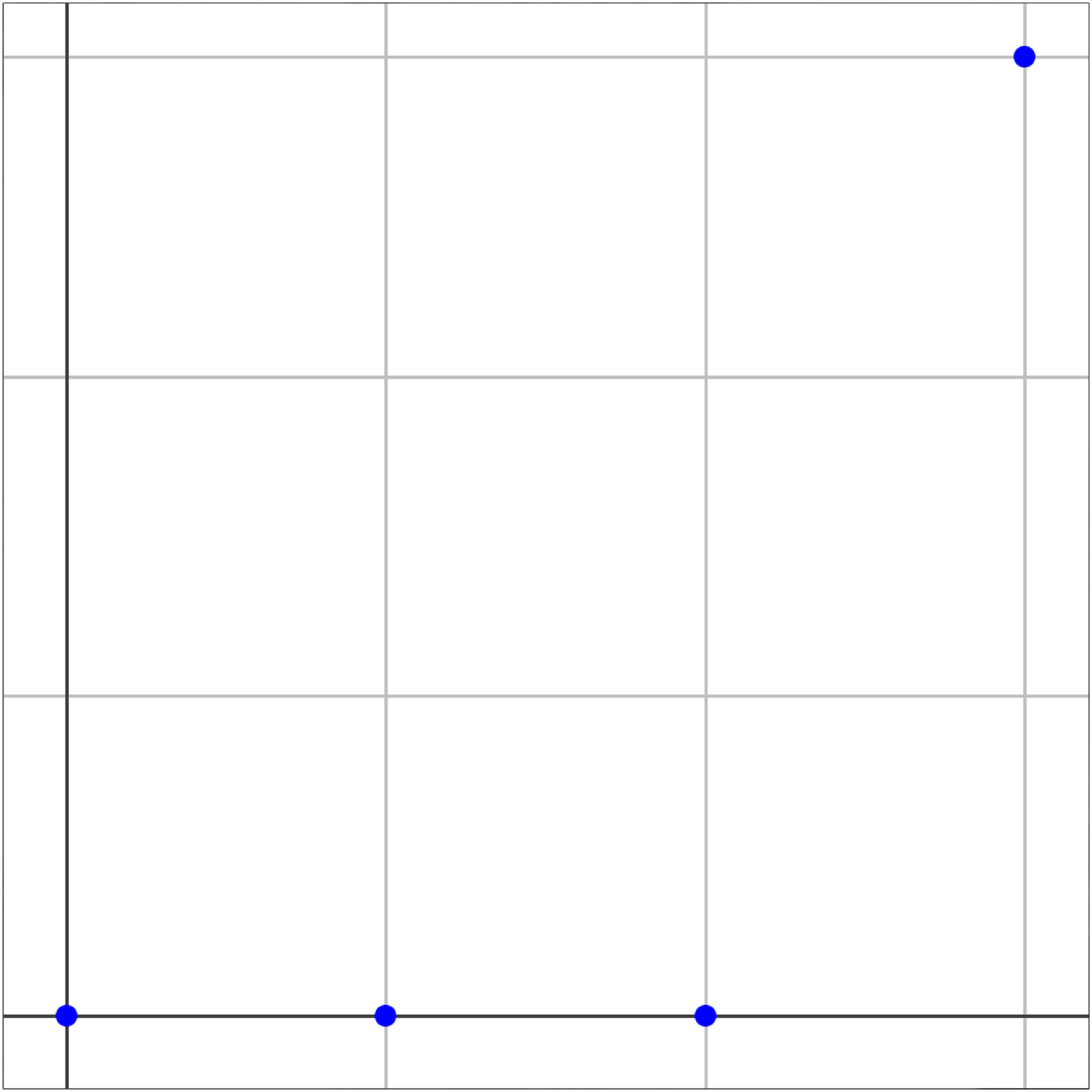} &
 \includegraphics[width=0.17\textwidth]{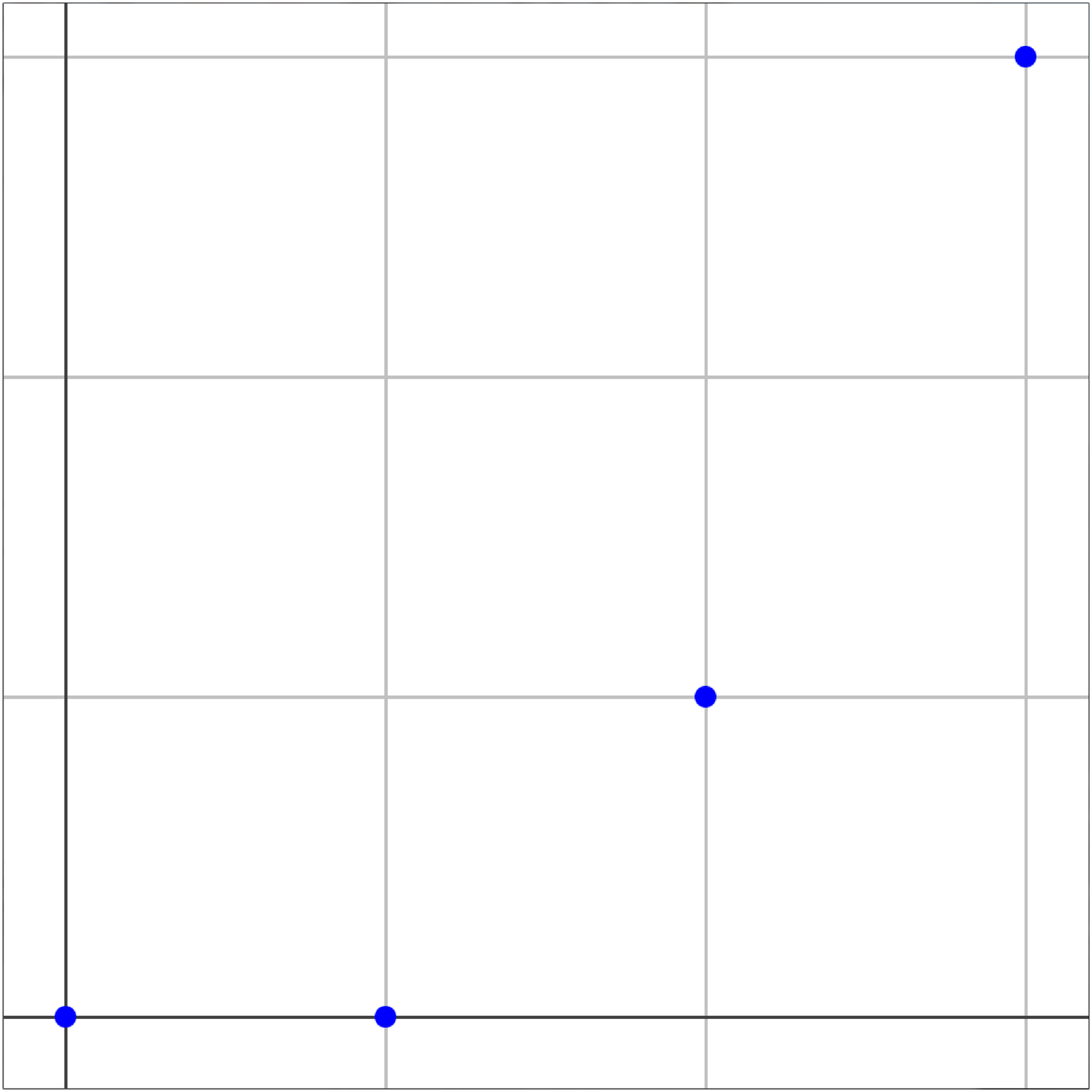} \\
 Point Set 3 & Point Set 4 \\
 \footnotesize{(0,0), (1,0), (2, 0), (3, 3)} &
 \footnotesize{(0,0), (1,0), (2, 1), (3, 3)}
 \end {tabular}
 \caption{The four-point sets used to compare elastic energies and curvature variations of various spline curves.}
 \label{fig:four-point-sets}
\end{figure}

We compare the elastic energies of spline segments between the center two points of various four-point sets in the plane as shown in Figure \ref{fig:four-point-sets}.

The elastic energies and curvature variations of our spline and standard splines with several parameter choices, as computed using $E = \int_{t_2}^{t_3} \kappa^2 \, \mbox{d}t,$ and $V = \int_{t_2}^{t_3} \dot{\kappa}^2 \, \mbox{d}t,$ are shown in Table \ref{tab:energies}.

\begin{table}[!ht]
 \centering
\begin{tabular}{l r r r r r r}
  & \rotatebox{90}{\textbf{This Paper's}} 
    \rotatebox{90}{\textbf{Method}} &
    \rotatebox{90}{\textbf{Catmull-Rom}} &
    \rotatebox{90}{\textbf{Cardinal}}
    \rotatebox{90}{($\tau=0.1$)} &
    \rotatebox{90}{\textbf{Cardinal}}
    \rotatebox{90}{($\tau=0.5$)} &
    \rotatebox{90}{\textbf{Kochanek-Bartels}}
    \rotatebox{90}{($\beta=0.5$)} &
    \rotatebox{90}{\textbf{Kochanek-Bartels}}
    \rotatebox{90}{($\beta=-0.5$)} \\
  \hline
  \textbf{Set 1} \\
  E: &  6.83 &  13.46 &  16.90 &   71.53 &  15.73 &  9.41 \\
  V: &   465 &   1742 &   2684 &   41012 &   2689 &   943 \\
  \hline
  \textbf{Set 2} \\
  E: &  0.50 &   0.50 &   0.66 &    3.65 &   0.85 &  0.85 \\
  V: &  25.6 &   25.6 &     55 &    2194 &   15.1 &  15.1 \\
  \hline
  \textbf{Set 3} \\
  E: &  0.49 &   4.00 &   3.69 &    4.93 &   1.57 &  6.62 \\
  V: &  14.6 &   81.6 &   78.6 &     404 &   16.4 &   212 \\
  \hline
  \textbf{Set 4} \\
  E: &  0.45 &   0.17 &   0.19 &    0.90 &   0.69 &  0.12 \\
  V: &   3.0 &    0.3 &    1.8 &     351 &   10.2 &   1.1 \\
  \hline
\end{tabular}
\caption{Energies and curvature variations of various splines, over the interval $t_2 < t < t_3$.}
\label{tab:energies}
\end{table}

\begin{figure}[!ht]
 \centering
 \begin{tabular}{c c}
  \includegraphics[width=0.17\textwidth]{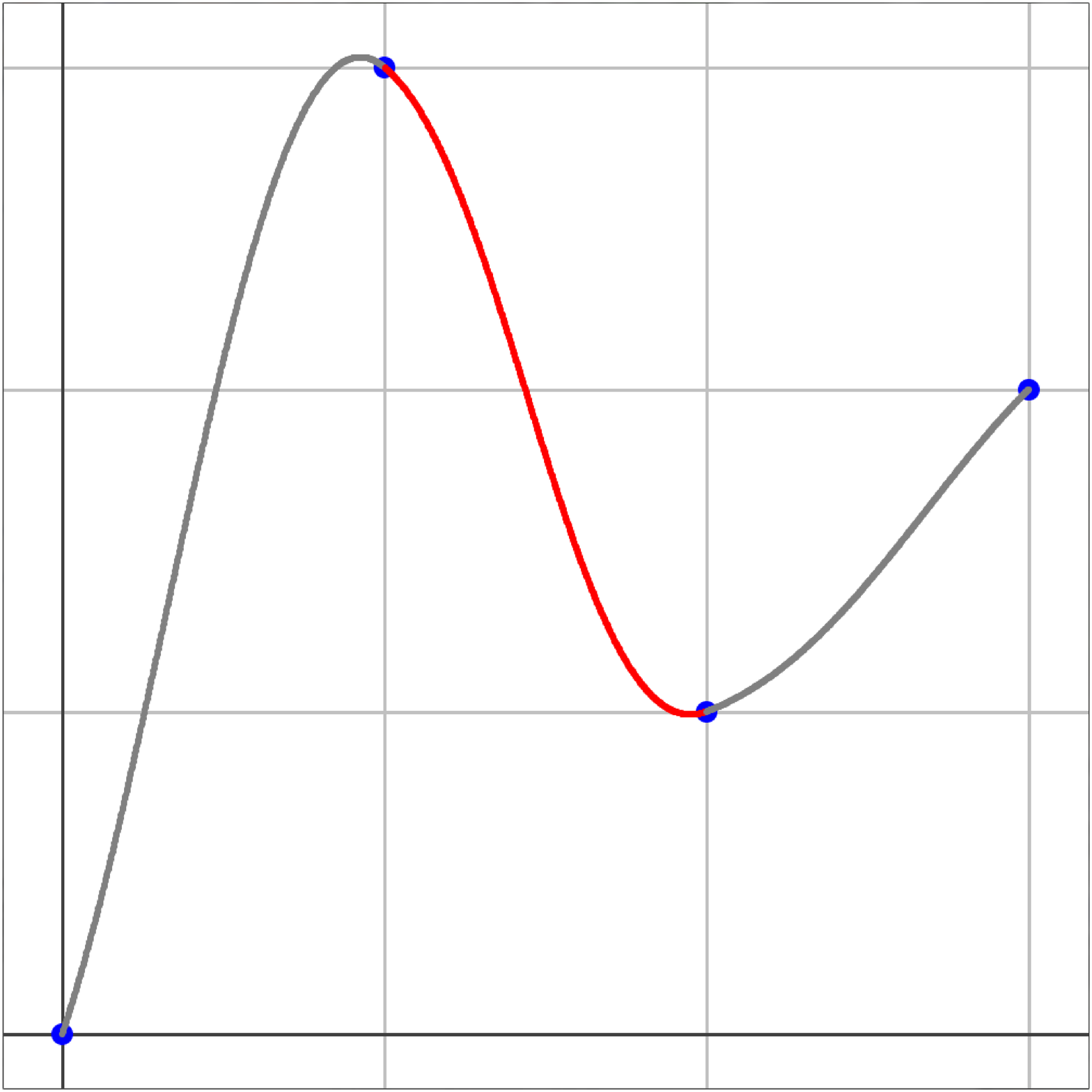} &
  \includegraphics[width=0.17\textwidth]{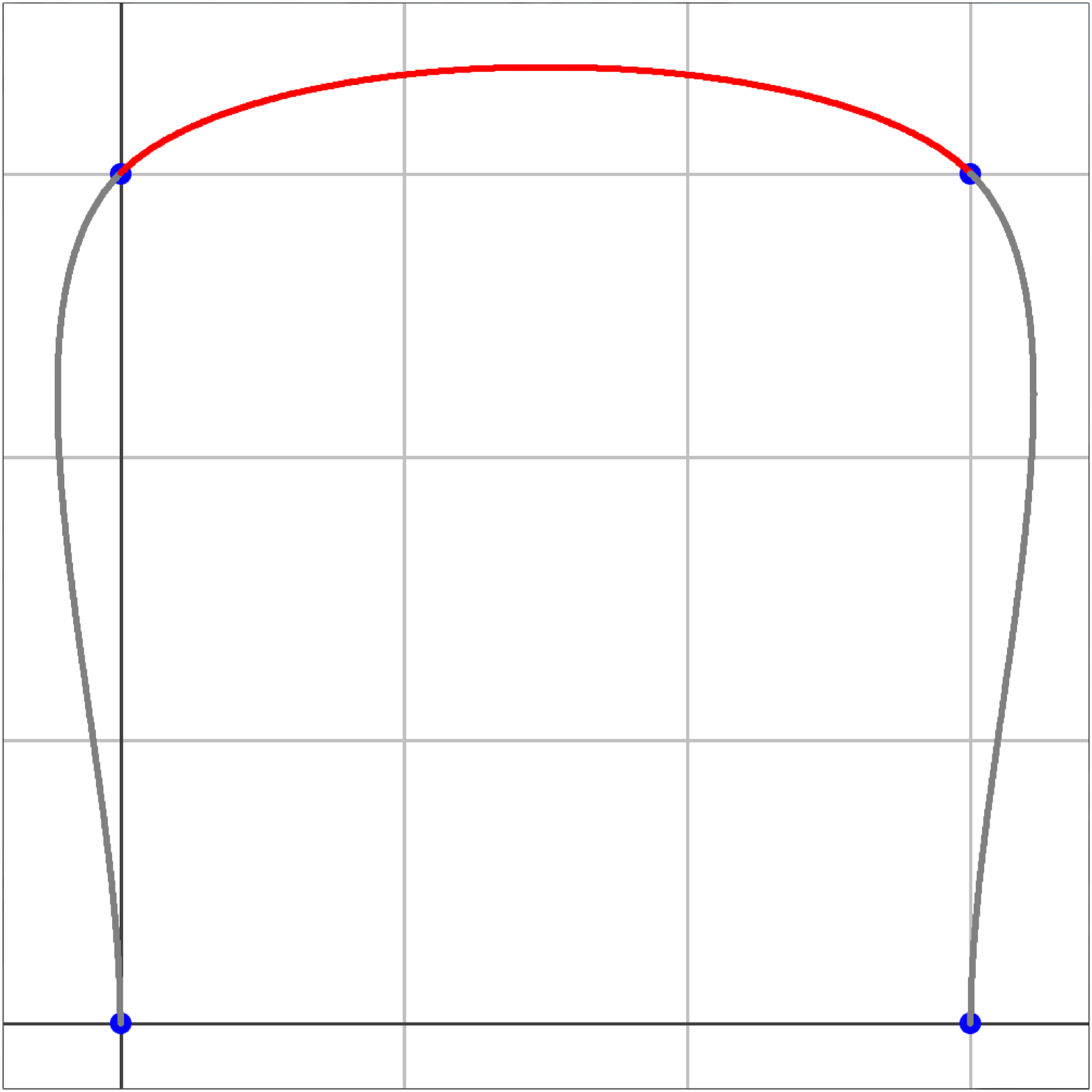} \\
  \includegraphics[width=0.17\textwidth]{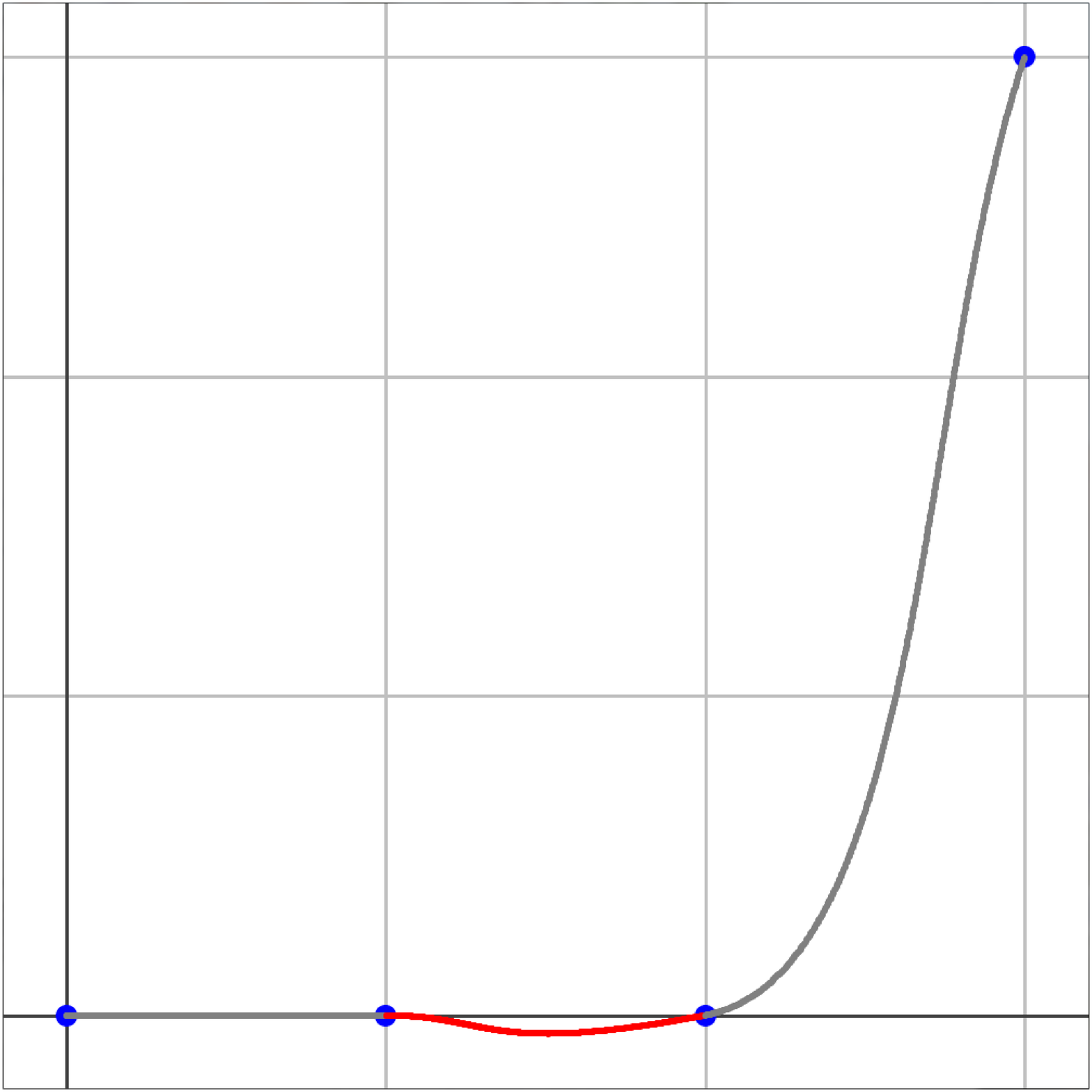} &
  \includegraphics[width=0.17\textwidth]{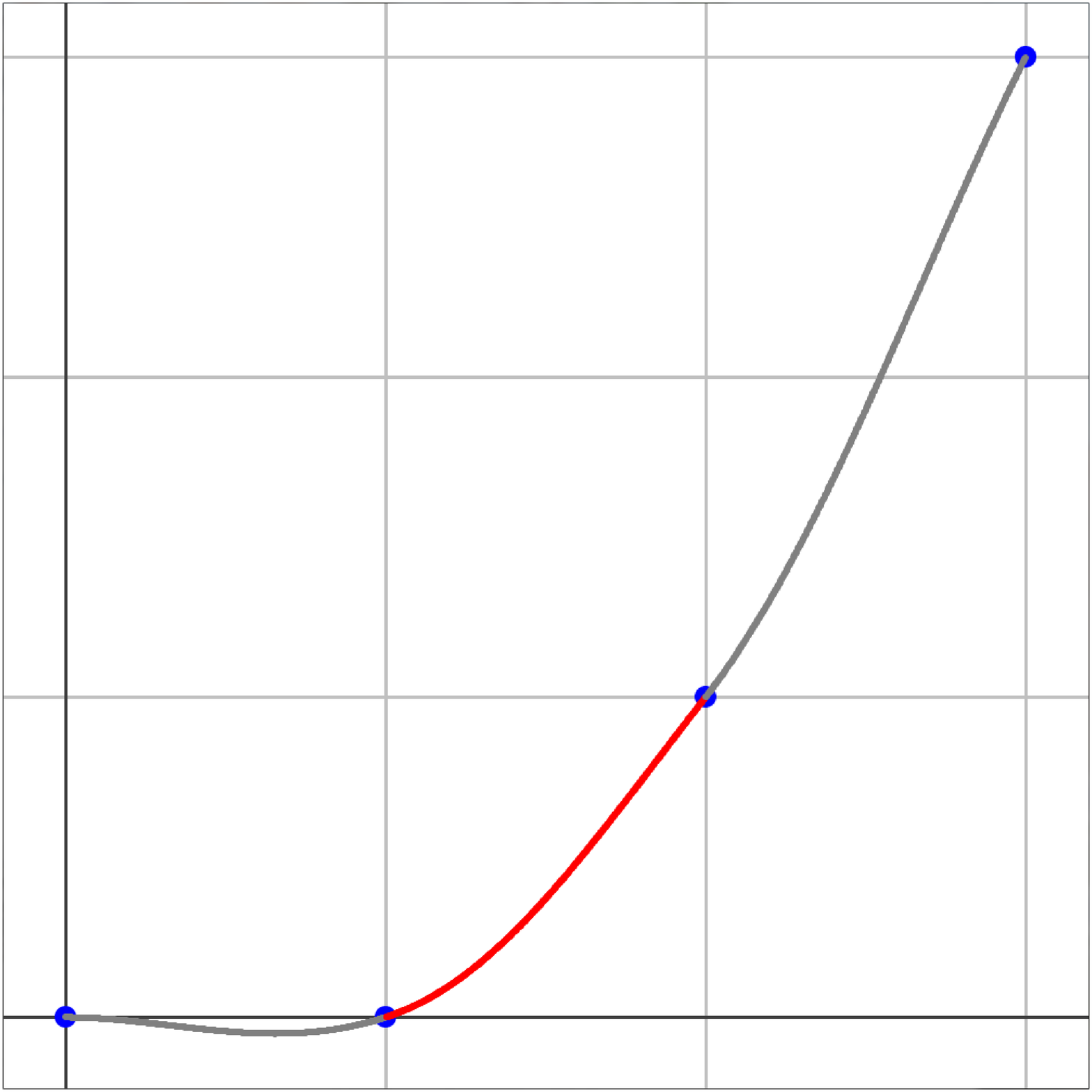}
 \end {tabular}
 \caption{The splines generated using the minimum-energy quadratic for the four test point sets. The energies in Table \protect{\ref{tab:energies}} are computed over the portion shown in red.}
 \label{fig:my-splines}
\end{figure}

\begin{figure}[!ht]
 \centering
 \includegraphics[width=0.3\textwidth]{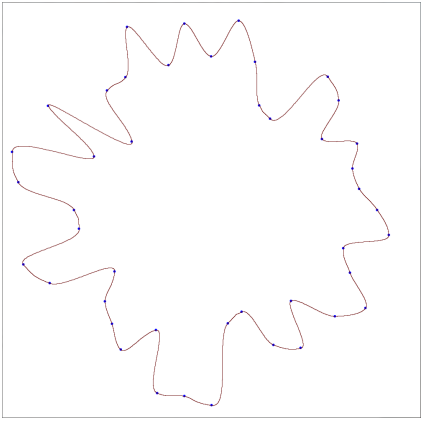}
 \caption{A randomly-generated ordered point set and the interpolating spline generated by building Hermite splines between successive points where the tangent vectors at each point are those of the least-energy quadratic through the point and its two neighbors.}
\label{fig:sample_point_set}
\end{figure}

The corresponding spline curves generated using min\-i\-mum-energy quadratics are shown in Figure \ref{fig:my-splines}.
An example of a complete point set with the spline generated by the method outlined above is shown in Figure \ref{fig:sample_point_set}.

\section{Conclusion}

The method we have demonstrated in this paper to generate Hermite splines demonstrates significantly lower energies and total curvature variations than other commonly used implementations for many (but not all) point sets.  The calculations involved are more complex, to be sure, but these calculations need only be done during des\-ign-time; calculating with or rendering the resulting spline curves can be accomplished using existing highly optimized Hermite curve routines.  In interactive applications, this method can also be used to provide a reasonable default set of tangent vectors at intermediate control points which act as a starting point for subsequent user manipulation.

We acknowledge that there are point sets for which the splines generated by this method have higher energy and curvature variation than those generated by other methods such as Catmull-Rom.  Having this method available, however, would allow a designer to compute the energies for each potential spline implementation and select the optimal spline based on energy and curvature variation, which in many cases would be the spline whose tangent vectors are found using our method.

\appendix
\section{Implementation algorithm} \label{appendix:implementation}

This appendix gives an algorithm for Hermite spline generation based on the minimum-energy quadratic.

\subsection*{Inputs:}
Four points $\bfp_i = (x_i, y_i)$, $i=1,2,3,4$, and a value of interpolation parameter $t$ in the range $[0,1]$,
where $t=0$ corresponds to point $\bfp_2$ and $t=1$ corresponds to point $\bfp_3$.

\subsection*{Output:}
The interpolated point $\bfr(t) = (x_t, y_t)$ on the segment of the curve between $\bfp_2$ and $\bfp_3$.

\subsection*{Algorithm:}

\noindent Procedure 1: To be performed when the set of control points changes.

\footnotesize
\noindent
\rule{\linewidth}{0.2mm}
\begin{alltt}
{\color{blu}FOR EACH} consecutive sequence of four points \(\mathbf{p}\sb{1},\mathbf{p}\sb{2},\mathbf{p}\sb{3},\mathbf{p}\sb{4}\)
    {\color{blu}SET} \(\mathbf{v}\sb{2}\) = COMPUTE_TAN(\(\mathbf{p}\sb{1},\mathbf{p}\sb{2},\mathbf{p}\sb{3}\))
    {\color{blu}SET} \(\mathbf{v}\sb{3}\) = COMPUTE_TAN(\(\mathbf{p}\sb{2},\mathbf{p}\sb{3},\mathbf{p}\sb{4}\))
{\color{blu}END FOR}

{\color{blu}FUNCTION} COMPUTE_TAN(\(\mathbf{q},\mathbf{r},\mathbf{s}\))
    {\color{grn}// Shift first point to origin, scale so \(|\tilde{\mathbf{s}}|=1\)}
    \(\tilde{\mathbf{r}}=\left(\mathbf{r}-\mathbf{a}\right)/|\mathbf{s}-\mathbf{a}|, \tilde{\mathbf{s}}=\left(\mathbf{s}-\mathbf{a}\right)/|\mathbf{s}-\mathbf{a}|\)

    {\color{grn}// Rotate such that point \(\tilde{\mathbf{s}}\) goes to (1,0)}    
    \(\widehat{\mathbf{q}}=\mathbf{0}\), \(\widehat{\mathbf{r}}=(\tilde{\mathbf{r}}\sb{x}\tilde{\mathbf{s}}\sb{x}+\tilde{\mathbf{r}}\sb{y}\tilde{\mathbf{s}}\sb{y}, -\tilde{\mathbf{r}}\sb{x}\tilde{\mathbf{s}}\sb{y}+\tilde{\mathbf{r}}\sb{y}\tilde{\mathbf{s}}\sb{x})\), \(\widehat{\mathbf{s}}=(1, 0)\)
    
    {\color{grn}// Compute the three roots of the cubic}
    \(\beta=1-2\widehat{\mathbf{r}}\sb{x}\), \(\gamma=(4(\widehat{\mathbf{r}}\sb{x}-|\widehat{\mathbf{r}}|\sp{2})-3)\sp{3}/27\), \(l=\sqrt[6]{-\gamma}\)
    \(\phi\sb{1}=\tan\sp{-1}\left(\sqrt{-\gamma-\beta\sp{2}}/\beta\right)/3\)
    \(\mu\sb{r}=l\,\cos\phi\sb{1}\), \(\mu\sb{i}=l\,\sin\phi\sb{1}\)
    \(\phi\sb{2}=\tan\sp{-1}\left(-\sqrt{-\gamma-\beta\sp{2}}/\beta\right)/3\)
    \(\zeta\sb{r}=l\,\cos\phi\sb{2}\), \(\zeta\sb{i}=l\,\sin\phi\sb{2}\)
    \(T\sb{1}=1/2+\mu\sb{r}+\zeta\sb{r}/2\)
    \(T\sb{2}=1/2-1/4[\mu\sb{r}+\zeta\sb{r}+\sqrt{3}(\mu\sb{i}-\zeta\sb{i})]\)
    \(T\sb{3}=1/2-1/4[\mu\sb{r}+\zeta\sb{r}-\sqrt{3}(\mu\sb{i}-\zeta\sb{i})]\)

    {\color{grn}// Find the root between 0 and 1}
    {\color{blu}IF} \(T\sb{1}>0\) AND \(T\sb{1}<1\) {\color{blu}THEN}
        \(T=T\sb{1}\)
    {\color{blu}ELSEIF} \(T\sb{2}>0\) AND \(T\sb{2}<1\) {\color{blu}THEN}
        \(T=T\sb{2}\)
    {\color{blu}ELSEIF} \(T\sb{3}>0\) AND \(T\sb{3}<1\) {\color{blu}THEN}
        \(T=T\sb{3}\)
    {\color{blu}END IF}
    
    {\color{grn}// Construct the \(\mathbf{a}\sb{1}\) and \(\mathbf{a}\sb{2}\) coefficients}
    \(\mathbf{a}\sb{1}=[(\mathbf{r}-\mathbf{q})-T(\mathbf{s}-\mathbf{q})]/(T\sp{2}-T)\)
    \(\mathbf{a}\sb{2}=\mathbf{s}-\mathbf{q}-\mathbf{a}\sb{1}\)
    
    {\color{grn}// Build and return the tangent vector \(\mathbf{v}\)}
    {\color{blu}RETURN} \(\mathbf{v}=2\mathbf{a}\sb{1}T+\mathbf{a}\sb{2}\)
\end{alltt}

\noindent Procedure 2: Hermite cubic spline interpolation, performed for each interpolation, with a value of \(t\) between 0 and 1.

\footnotesize
\noindent
\rule{\linewidth}{0.2mm}
\begin{alltt}
\(\mathbf{r}(t)=(2t\sp{3}-3t\sp{2}+1)\mathbf{p}\sb{2}+(-2t\sp{3}+3t\sp{2})\mathbf{p}\sb{3}+(t\sp{3}-2t\sp{2}+t)\mathbf{v}\sb{2}+(t\sp{3}-t\sp{2})\mathbf{v}\sb{3}\)
\end{alltt}
\rule{\linewidth}{0.2mm}

\bibliographystyle{elsarticle-num}
\bibliography{MinEnergyQuadratic}

\end{document}